\documentclass[12pt,oneside]{amsart}
\usepackage{graphicx, epsfig}

\newfont{\sheaf}{eusm10 scaled\magstep1}

\newcommand{\ra}{\ensuremath{\rightarrow}}

\def\eea{\end{eqnarray*}}
\def\bea{\begin{eqnarray*}}

\newcommand{\Proof}{{\it Proof. }}

\newtheorem{teo}{Theorem}[section]
\newtheorem{df}[teo]{Definition}
\newtheorem{lem}[teo]{Lemma}
\newtheorem{cor}[teo]{Corollary}
\newtheorem{ex}[teo]{Example}
\newtheorem{oss}[teo]{Remark}
\newtheorem{prop}[teo]{Proposition}

\newcommand{\suchthat}{\mathrel{|}}

\newcommand{\C}{\ensuremath{\mathbb{C}}}
\newcommand{\R}{\ensuremath{\mathbb{R}}}

\newcommand{\F}{\ensuremath{\mathbb{F}}}

\newcommand{\hol}{\ensuremath{\mathcal{O}}}

\newcommand{\PP}{\ensuremath{\mathbb{P}}}

\newcommand{\PPP}{\ensuremath{\mathcal{P}}}

\newcommand{\sE}{{\mathcal E}}

\begin{document}

\title[Real singular Del Pezzo surfaces]{Real singular Del Pezzo surfaces and threefolds
fibred by  rational curves, I.}

\author{Fabrizio Catanese\\
  Universit\"at Bayreuth \\ Fr\'ed\'eric Mangolte\\ Universit\'e de Savoie }

\thanks{The research of the  authors was  supported by the SCHWERPUNKT "Globale Methoden in
der komplexen Geometrie" and the ANR grant "JCLAMA" of the french "Agence Nationale de la Recherche".}

\date{july 24, 2007}

\begin{abstract}
Let $W \to X$ be a real smooth projective threefold fibred by rational curves. Koll\'ar proved that if $W(\R)$ is orientable a connected component $N$ of $W(\R)$ is essentially  either a Seifert fibred  manifold or a connected sum of lens spaces.
Let $k : = k(N)$ be the integer defined as follows:
If $ g \colon N \to F$ is a Seifert fibration, one defines $k : = k(N)$ as the number of  multiple fibres of $g$, while, if $N$ is a connected sum of lens spaces, $k$ is defined as the number of lens  spaces different from $\mathbb{P}^3(\mathbb{R})$. Our
Main Theorem says:
If $X$ is a geometrically rational surface, then $k \leq 4$.
Moreover we show that if $F$ is diffeomorphic to $S^1 \times S^1$, then
 $W(\R)$ is connected and  $k = 0$.

These results answer in the affirmative two questions  of Koll\'ar who proved  in 1999 that $k \leq 6$ and suggested that 4 would be the sharp bound.
We derive the Theorem from a careful study of real singular Del Pezzo surfaces with only Du Val singularities.
\end{abstract}

\maketitle

\begin{quote}\small
\end{quote}

\tableofcontents

\bigskip

%\newpage

\section*{INTRODUCTION}

In complex algebraic geometry there is an established principle that the
Kodaira dimension of a smooth complex projective variety $W$ of dimension $n$ strongly
influences the topology of the set $ W(\C)$ of its complex points.
This principle is clearly manifest already in dimension 1, and related
to  other points of view, as the uniformization theorem, and the concept of curvature.
This principle, although in a more difficult and complicated way, still goes on to hold
in higher dimensions.

Let us assume now that $W$ is a smooth real projective variety and let us consider
the topology of the set $ W(\R)$ of its real points. In dimension 1, the connected
components are just diffeomorphic to the circle $S^1$, and their number is not dictated
by the genus (there is only the Harnack inequality which gives $g+1$ as upper bound
for the number of connected components of  $ W(\R)$).

So, there had been for some time the belief that the Kodaira dimension of $W$
would not affect at all the topology of a connected component $N$ of  $ W(\R)$.
This belief is contradicted already by the example of real algebraic surfaces
of nonpositive Kodaira dimension (see for instance \cite{Co14}, \cite{Si89}, \cite{dik}
and \cite{kolIP}).

In a very interesting series of papers (\cite{KoI,KoII,KoIII,KoIV}) J\'anos Koll\'ar used the
recent progress on the minimal model program for threefolds in order to understand the
topology of the connected components $  N \subset   W(\R)$, especially in the case
where $W$ has Kodaira dimension  $ - \infty$.

Our note takes the origin from some questions that Koll\'ar set in the third
article of the series (\cite{KoIII}), and we prove some optimal estimates
that Koll\'ar conjectured to hold.

The present note is mainly devoted to the proof of the following 

\begin{teo}\label{teo:main} Let $X$ be a projective surface defined over $\R$. Suppose that
$X$ is geometrically rational with Du Val singularities. Then a connected component $M$ of
the topological normalization $\overline {X(\R)}$ contains at most 4 Du Val singular points
which are either not of type $A^-$ or of type $A^-$ but globally separating.
\end{teo}

Applying this  result to rational curve fibrations over rational surfaces, we obtain the
answer to two of the three questions set by Koll\'ar (remark 1.2 of \cite{KoIII}).
In a second note,  with slightly different methods, we plan to answer also the third
question.

Let us now explain these applications in more detail.

Let $f \colon W \to X$ be a real smooth projective threefold fibred  by rational curves.
Suppose that $W(\mathbb{R})$ is orientable. Then, by 
\cite[Theorem~1.1]{KoIII}, a connected component $N \subset W(\mathbb{R})$ is a Seifert
fibred manifold, or a  connected sum of lens spaces, or obtained from one of the above by
taking connected sums with a finite number of copies of $\mathbb{P}^3(\mathbb{R})$ and a
finite  number of copies of $S^1\times S^2$. Note that 
 in \cite{hm1} and \cite{HM05b} it was shown that all the manifolds $N$ as above
do indeed occur.

Note also that the connected sum $N_1 \# N_2$ is
taken in the category of oriented manifolds, where in general $N_1 \# N_2$ is not
homeomorphic to $N_1 \# -N_2$. But for the particular choice $N_2 =
\mathbb{P}^3(\mathbb{R})$ or $N_2 = S^1\times S^2$, the connected sums $N_1 \# N_2$ and $N_1
\# -N_2$ are diffeomorphic, see e.g. \cite{Hem76}.

Take a decomposition $N = N' \#^a\mathbb{P}^3(\mathbb{R})\#^b(S^1\times S^2)$ with $a + b$
maximal and observe that this decomposition is unique by a theorem of Milnor \cite{Mil62}.

We shall focus our attention on  the integer $k : = k (N)$   defined as follows:

\begin{enumerate}
\item if  $g \colon N' \to F$ is a Seifert fibration,  $k$ denotes  the number of multiple
fibres of
  $g$;
\item if $N'$ is a connected sum of lens spaces,  $k$ denotes the number of lens spaces.
\end{enumerate}

Observe that when $N'$ is a connected sum of lens spaces, the number $k$ is well defined
(again by Milnor's theorem).

We can then apply the result of Theorem \ref{teo:main} concerning singular rational surfaces
in order to  answer one question of Koll\'ar,  \cite[Remark 1.2 (1)]{KoIII}.

\begin{cor}\label{cor:multiple} Let $ W \to X$ be a real smooth projective $3$-fold fibred
by rational curves over a geometrically rational surface $X$. Suppose that $W(\R)$ is
orientable. Then for each connected component $N \subset W(\R)$, $k(N)  \leq 4$.
\end{cor}

 Note that Koll\'ar showed in \cite{KoIII} the optimality of the above estimate
in case 1).

The following theorem answers another question of Koll\'ar,  \cite[Remark 1.2
(3)]{KoIII}

\begin{teo}\label{teo:torus} Let $W$ be a real smooth projective $3$-fold fibred by rational
curves over a geometrically rational surface $X$. Suppose that the fibration is defined over
$\R$ and that $W(\R)$ is orientable. Let $N \subset W(\R)$ be a connected component which
admits a Seifert fibration $g \colon N \to S^1\times S^1$. Then $g$ has no multiple fibres.
Furthermore, $X$ is then rational over $\R$ and $W(\R)$ is connected.
\end{teo}

Section 1 is devoted to recalling the basic notions which come into play, especially
the local and global separation properties of Du Val singularities of types
$A_{\mu}$. Hence the basic invariants $m_i$ of a Du Val surface are defined.

Section 2 is the heart of the paper and contains a detailed description of the
topological normalization of a  real Du Val Del Pezzo surface with more than
4 singular points. The description is based on the classical representation 
of the quadric cone on the plane which transforms the hyperplane sections of
the cone to  parabolae whose axis has a given direction.

Section 3 proves Corollary \ref{cor:multiple}, while Section 4 is devoted 
to the proof of Theorem \ref{teo:torus}.

%%%%%%%%%%%%%%%%%%%%%%%%%%%%%%%%%%%%%%%%
\section{Real Du Val surfaces}\label{sec:dvsurf}
%%%%%%%%%%%%%%%%%%%%%%%%%%%%%%%%%%%%%%%%

The aim of this section is to reduce the proof of Theorem~\ref{teo:main} to the study of a
certain kind of rational surfaces. The first part is close to the treatment in
\cite[Section~9]{KoIII}.

On a surface, a rational double point is called a Du Val singularity.  Over $\mathbb{C}$,
these singularities are classified by their Dynkin diagrams, namely 
$A_\mu$, $\mu \geq 1$, $D_\mu$, $\mu \geq 4$, $E_6$,
$E_7$, $E_8$.

Over $\mathbb{R}$, there are more possibilities. In particular, a surface singularity will
be said to be {\em of type $A^+_\mu$} if it is real analytically equivalent to
$$ x^2 + y^2 - z^{\mu + 1} = 0,\ \mu \geq 1\;;
$$ and {\em of type $A^-_\mu$} if it is real analytically equivalent to
$$ x^2 - y^2 - z^{\mu + 1} = 0,\ \mu \geq 1\;.
$$ The type $A_1^+$ is real analytically isomorphic to $A_1^-$; otherwise, singularities
with different names are not isomorphic. For $\mu$ odd, there is another real form of
$A_\mu$ given by $x^2 + y^2 + z^{\mu + 1} = 0$. We exclude this type of singular point
because an isolated real point gives rise to $\emptyset$ on the minimal resolution.

\begin{df} Let $X$ be a projective surface. The surface $X$ is called a {\em Du Val} surface
if  $X$ has only rational double points as singularities.
\end{df}

We want to use a suitable minimal model for $X$. In the minimal model program for real Du
Val surfaces, the most useful statement for our purpose is the following description of the
extremal contractions.

\begin{teo}\cite[Th.~9.6]{KoIII}\label{teo:extremal ray} Let $X$ be a real Du Val surface, 
$\overline{NE}(X)$ be its cone of curves, and $R \subset \overline{NE}(X)$ be a
$K_X$-negative extremal ray. Then $R$ can be contracted. Furthermore, if $c \colon X \to Y$
is the contraction, $c$ is one of the following:
\begin{itemize}
\item $Y$ is a Du Val surface, $c$ is birational, and $\rho(Y) = \rho(X) - 1$,
\item $Y$ is a smooth curve, $\rho(X) = 2$, and $c \colon X \to Y$ is a conic bundle,
\item $Y$ is a point, $\rho(X) = 1$ and $X$ is a Du Val Del Pezzo surface (i.e. $-K_X$ is
ample). 
\end{itemize}

\end{teo}

To apply the minimal model program for our purposes, we need to understand the behavior of
$c$ when $c$ is a birational contraction. We begin with a typical example.

\begin{ex}\label{ex:weightedblowup} Let $Y$ be a real Du Val surface, $x \in Y$ be a smooth
real point, and $\mu > 0$ be an integer. Blow-up $Y$ at $x$, and denote by $E_0$ the
exceptional curve of the blow-up $Y_0 \to Y$. Then take repeatedly the blow-up $Y_{l+1} \to
Y_l$ at a general point on the exceptional curve $E_{l}$ for $l =  0,1,\dots,\mu - 1$. The
exceptional divisor of the composition of blow-ups $Y_\mu \to Y$ is a chain of rational
curves whose configuration is of the form:

 \begin{figure}[htbp]
\begin{center}
   \epsfysize=0.6cm
  \epsfbox{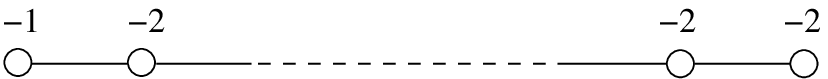}
\end{center}
\end{figure}

Contracting the $(-2)$-curves $E_{\mu-l},\ l =  0,1,\dots,\mu - 1$, we get a surface $X$
with a singularity of type $A_\mu^-$.
\end{ex}

The interesting fact is that the birational contractions of Theorem~\ref{teo:extremal ray}
involve only this kind of construction (see \cite[Th.~9.6]{KoIII}).

As we shall see in Section~\ref{sec:3fold}, bounding the number of certain singularities on $X(\R)$ yields
a bigger upper bound for $k(N)$ than the one stated in Corollary~\ref{cor:multiple}.
In order to obtain this finer estimate we have to bound this number on each component of the topological normalization $\overline{X(\mathbb{R})}$ of $X(\R)$, which we are going now to define.

\begin{df} Let $V$ be a simplicial complex with only a finite number of points 
$x \in V$ where $V$ is not a manifold. Define the {\em topological normalization}
$$
\overline{n} \colon \overline{V} \to V
$$ as the unique proper continuous map such that $\overline{n}$ is a  homeomorphism over the
set of points
  where $V$ is a manifold and $\overline{n}^{-1}(x)$ is in one-to-one  correspondence with
the connected components of a good punctured neighborhood of $x$ in $V$ otherwise.
\end{df}

Observe that if $V$ is pure of dimension 2, then $\overline{V}$ is a topological manifold.
Indeed each point of $\overline{V}$ has a neighbourhood which is a cone over $S^1$.

\begin{df} Let $X$ be a real algebraic surface with isolated singularities, and let
$x \in X(\mathbb{R})$ be a singular point of type $A^\pm_\mu$ with $\mu$ odd. The
topological normalization 
$\overline{X(\mathbb{R})}$ has two connected components locally near $x$. We will say that
$x$ is {\em globally separating} if  these two local components are on different connected
components of $\overline{X(\mathbb{R})}$ and {\em globally  nonseparating} otherwise.
\end{df}

One can produce an arbitrarily high number of singular points of type $A_\mu^-$ by the
construction of Example~\ref{ex:weightedblowup}, but these singular points are globally
nonseparating. Indeed, when $\mu$ is even, the singular point is in fact  locally
nonseparating, and when $\mu$ is odd, then the inverse image of the last $S^1 = E_\mu(\R)$
yields a segment in $\overline{X(\mathbb{R})}$ connecting the two points. The key point for
the sequel is the next lemma.

\begin{df} Let $X$ be a real Du Val surface, let
\begin{multline}
\PPP_X := \operatorname{Sing}X \setminus \left\{ x \textrm{ of type } A^-_\mu,\ \mu   \textrm{ even}\right\} \\
\setminus \left\{ x \textrm{ of type } A^-_\mu,\ \mu \textrm{ odd and
} x \textrm{ is globally nonseparating}\right\} \;.  
\end{multline}
\end{df}

We have 

\begin{lem}\cite[Cor.~9.7]{KoIII}\label{lem:m_i} Let $X$ be a real Du Val surface, let
$\overline{n} \colon \overline{X(\R)} \to X(\R)$ be the topological normalization, and let
$M_1,M_2,\dots, M_r$ be the connected components of  $\overline{X(\R)}$. The unordered
sequence of numbers $ m_i := \#(\overline{n}^{-1}(\PPP_X) \cap M_i )$, $i = 1,2,\dots, r$ is
an invariant of extremal birational contractions of Du Val surfaces.
\end{lem}

By Theorem~\ref{teo:extremal ray} and Lemma~\ref{lem:m_i}, it suffices to prove
Theorem~\ref{teo:main} in the case when $X$ is a conic bundle or a Del Pezzo surface with
$\rho(X) = 1$. Conic bundles were analysed in  \cite[Section~9]{KoIII}. 
The remaining case is when $X$ is a Del Pezzo surface. We
still slightly reduce the problem to the case where $X$ is a degree 1 Del Pezzo surface.

\begin{lem}\label{lem:reducetodps1} Let $X$ be a real Du Val Del Pezzo surface possessing a smooth
real point and having $\rho(X) = 1$. Then there exists a blow-up of $X$ in smooth points
yielding $Z$ which is a conic bundle if $\deg X \geq 3$. Else we get $Z$ a singular Del
Pezzo surface of degree 1 with $\rho(Z) \leq 2$.   
\end{lem}

\begin{proof} Set $d := \deg X$. If $d \geq 3$, blow-up $(d - 3)$ smooth points until you
get a real cubic surface $Z$. The surface $Z$ contains a real rational line $L$. We get  $L \subset Z \subset \PP^3$, and
$\pi_L \colon \PP^3 - L \to \PP^2$   is a morphism and yields a real conic bundle.
If $d = 2$ blow-up a smooth real point : $\rho(X)$ increases by 1.
\end{proof}

%%%%%%%%%%%%%%%%%%%%%%%%%%%%%%%%%%%%%%%%
\section{Singular Del Pezzo surfaces of degree 1}\label{sec:sdps}
%%%%%%%%%%%%%%%%%%%%%%%%%%%%%%%%%%%%%%%%

Recall that  a Del Pezzo surface $X$ is by definition a  surface  whose anticanonical
divisor is ample.  We add the adjective {\em Du Val} to emphasize that we allow $X$ to have
Du Val  singularities (observe that for a Du Val surface, the canonical divisor is a Cartier
divisor).  Let $X$ be a real Du Val Del Pezzo surface  and let
$p \colon S \to X$ be the minimal resolution of singularities.  The smooth surface $S$ has
nef  anticanonical divisor $- K_S = p^*( -K_X)$, and is called a {\em weak Del Pezzo
surface}  by many authors. As we saw in Section~\ref{sec:dvsurf}, we can assume the Del
Pezzo surface $X$ to have degree 1 by blowing up a finite  number of pairs of conjugate
imaginary smooth points and some real smooth point (there are several choices to do this),
see Lemma~\ref{lem:reducetodps1}.  The anticanonical model of a Del Pezzo surface $X$ of
degree 1 is a ramified double covering $q \colon X \to Q$ of a quadric cone $Q \subset
\PP^3$ whose branch locus is the union of the vertex of the cone and a cubic section $B$ not
passing through the vertex, see e.g. \cite[Expos\'e V]{Dem80}.
 
Remark that the pull-back by $q$ of the vertex of the cone is a smooth point of $X$ and let
$X'$ be the singular elliptic surface obtained from $X$ by blowing up this smooth point. We
denote by $\overline{n} \colon \overline {X'(\R)} \to X'(\R)$ the topological normalization
of the real part.

% and by $S' \to X'$ the minimal resolution of singularities. 

We shall now make a series of considerations which will later lead to a proof of the
following.

\begin{prop}\label{prop:dps1case} For each connected component $M \subset \overline
{X'(\R)}$, 
$$
\#(\overline{n}^{-1}(\PPP_{X'}) \cap M) \leq 4 \;.
$$
\end{prop}

\bigskip

Recall that Hirzebruch surfaces are the $\PP^1$-bundles over $\PP^1$. The surface $X'$ is a
ramified double covering of the Hirzebruch surface $\F_2$ whose branch curve is the union of
the unique section of negative selfintersection, the section at infinity $\Sigma_\infty$,
and a trisection $B$ of the ruling $\F_2 \to \PP^1$ which is disjoint from $\Sigma_\infty$.  

The cone $Q$ is the weighted projective plane $\PP(1,1,2)$ with coordinates $(x_0,x_1,y_2)$,
and $X$ is the hypersurface in $\PP(1,1,2,3)$ with coordinates $(x_0,x_1,y_2,z)$ defined by
\begin{equation*} 
z^2 = y_2^3 + p_4(x_0,x_1)y_2 + q_6(x_0,x_1) \;.
\end{equation*}

We want to explain here the plane model of $Q$, in which the hyperplane sections of $Q$
embedded in $\PP^3$ by $H^0(\hol_Q(2))$ correspond to parabolae tangent to the line at
infinity $L_\infty = \{ w = 0\}$ at  the point $O :=  \{w = x = 0 \}$ of the projective
plane with coordinates $(x,y,w)$.  In other words, blow-up $O$ and then the infinitely near
point $O'$ to $O$ corresponding to the tangent of the line at infinity $L_\infty$, and
denote by $\tilde Q$ the resulting surface.  Denote by $E, E'$ the respective total
transforms of $O, O'$, and observe that $E = E' + E''$, $E''$ being a $(-2)$-curve.  The
linear system $H^0(\hol_{\tilde Q}(2H - E - E'))$ maps $\tilde Q$ birationally onto the
quadric cone $Q \subset \PP^3$, contracting the proper transform $\tilde L_\infty$ of
$L_\infty$ and $E''$ to points.  Since $\tilde L_\infty$ and $E''$ do not meet, first
contracting $\tilde L_\infty$ yields the Hirzebruch surface $\F_2$, whose $(-2)$-section
$\Sigma_\infty$ is the image of the curve $E''$.  Let us write everything using the
coordinates $(x,y,w)$ in $\PP^2$: then $H^0(\hol_Q(1))$ corresponds to $H^0(\hol_{\tilde
Q}(H - E))$ spanned by $w,x$, whereas $y_2 =: yw$ completes $w^2,wx,x^2$ to a basis of
$H^0(\hol_{\tilde Q}(2H - E - E')) \cong H^0(\hol_Q(2))$. Thus the morphism of $\tilde Q$ to
$\PP(1,1,2)$ is given by $x_0 := w$, $x_1 := x$, $y_2 := yw$.

The elliptic surface $X'$ is the double cover of $\F_2$ branched on $\Sigma_\infty$ and on
the curve $B$ corresponding to the curve of $Q$ of equation $y^3 + p_4(x_0,x_1)y +
q_6(x_0,x_1) = 0$. Thus the curve $B$ corresponds to the plane curve $w^3y^3 + p_4(w,x)yw +
q_6(w,x) = 0$ whose affine part has equation:
\begin{equation}\label{eq:B} y^3 + p_4(1,x)y + q_6(1,x) = 0 \;.
\end{equation}

Note that any parabola as above, i.e., a curve $C \in (2H - E - E')$ is disjoint from $E''$
(mapping to the vertex of the cone) unless it splits into two lines through the point $O$.
In particular, we may always change coordinates in the affine plane so that $C$ is
transformed into the line $y = 0$. In order to understand with coordinates the geometry at
infinity of parabolae as above, let us observe that $\F_2$ has two open sets isomorphic to
$\C \times \PP^1$. They have respective coordinates $\frac xw \in \C$, $(w,y) \in \PP^1$,
while on the other chart we have $\frac wx \in \C$, and homogeneous coordinates $(\frac
{x^2} w , y)$ (in fact $\frac {x^2} w /w = (\frac xw)^2$). The section $\Sigma_\infty$ at
infinity corresponds to the curve $E'' \subset \tilde X$ and is defined by $w = 0$ and
$\frac {x^2} w = 0$ on the respective charts. Then a parabola $yw = a_0w^2 + a_1 xw + a_2
x^2$ is given by the equation
$$
\frac 1\eta = a_0 + a_1 \frac xw + a_2(\frac xw)^2
$$ on the affine chart with coordinates $(\frac xw, \frac wy := \eta)$. Using these
coordinates at infinity it will be easy to see when some regions in the plane "meet" at
infinity in $\F_2$.

We shall now look for normal forms of Equation~ (\ref{eq:B}). Singular points of $X'(\R)$
are in one-to-one correspondence with singular points of $B(\R)$. The different cases we
shall now consider are distinguished by the number of irreducible components of the
trisection $B$.

%\medskip
\subsection*{Three components}

In this case, we shall see that any connected component of the topological normalization of
any real double cover ramified over $B$ will have at most 4 singular points. Observe that at least one of the three components is real.
Equation~(\ref{eq:B}) becomes 
$$ (y - \alpha(x))(y - \beta(x))(y - \gamma(x)) = 0
$$ and, changing real coordinates for $Q = \PP(1,1,2)$, we may assume $\gamma = 0$. The case $\beta = \overline{\alpha}$ where two components are complex conjugate leads to at most 2 singular points: $\operatorname{Re} \alpha(x) = 0, y = \operatorname{Im} \alpha(x)$. We can therefore assume that the three components are real.
Thus
Equation~(\ref{eq:B}) becomes 
$ (y - \alpha)(y - \beta)y  = 0
$ where $\alpha(x) = \alpha_0 + \alpha_1 x + \alpha_2 x^2$ and $\beta(x) = \beta_0 + \beta_1
x + \beta_2 x^2$ are polynomials of degree~2.

\begin{itemize}
\item Assume no 2 parabolae are tangent. Then, since we can permute the 3 curves, we can fix
the one which is the lowest at infinity (i.e., if we write the curves as $y = a_0 + a_1 x + a_2
x^2$, the one with the smallest value of $a_2$). Changing coordinates, we get only the curve $y = 0$
and two convex parabolae, i.e. with $\alpha_2 >0$ and $\beta_2 >0$.

The 6 intersection points are distinct and given by
$$ y = \alpha(x)\beta(x) = 0\;,\qquad
\alpha(x) = \beta(x) = y \;.
$$ The curve $B$ is real, thus if one of these points is not real, then the number of real
singular points is bounded above by 4 and we are done. From now on, we suppose that the six
points are real. Set 
\begin{equation}\label{eq:poly}
\begin{cases}
\alpha(x) = \alpha_2(x - a_1)(x - a_2),\ a_1 < a_2 \;;\\
\beta(x) = \beta_2(x - b_1)(x - b_2) \;.
\end{cases}
\end{equation}

Multiplying $y$ possibly by $\beta_2$, we may assume $\beta_2 = 1$.  We may reduce to the
case $0 < \alpha_2 < 1$ by possibly exchanging the roles of $\alpha$ and $\beta$. We can
further use a translation in the $x$ axis and assume $b_1 = - b_2$, then (\ref{eq:poly})
becomes:
$$
\begin{cases}
\alpha(x) = \alpha_2(x - a_1)(x - a_2),\ a_1 < a_2,\ 0 < \alpha_2 < 1 \;;\\
\beta(x) = (x^2 - b^2),\ 0 < b\;.
\end{cases}
$$

Up to reflection $x \leftrightarrow -x$, this leads to 4 possibilities, namely (see
Figure~\ref{fig:6points})
$$ b < a_1\;,\ -b < a_1 < b < a_2\;,\ a_1 < -b < b < a_2\;,\ -b < a_1 < a_2 < b \;.
$$ 

\begin{figure}[htbp]
\begin{center}
   \epsfysize=5.5cm
  \epsfbox{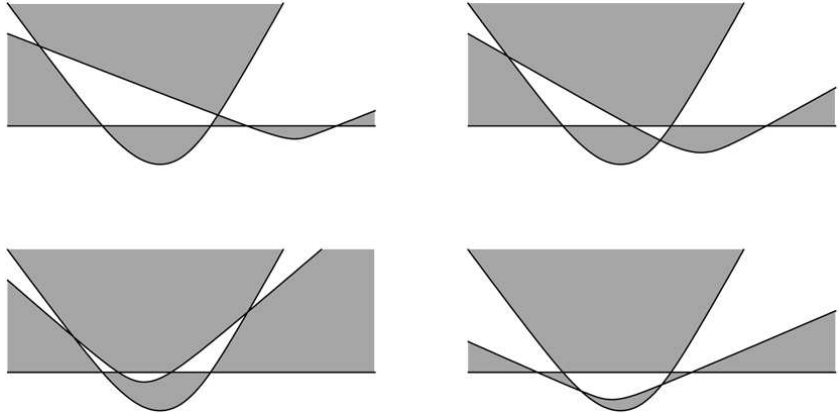}
\caption{Three parabolae, 6 singular points.}\label{fig:6points}
\end{center}
\end{figure}

\begin{oss}
Observe that in these figures two components are connected at infinity if their boundaries have two unbounded arcs belonging to the same pair of parabolae. 
\end{oss}

\medskip
\item Assume 2 parabolae are tangent. Then we cannot arbitrarily permute the 3 curves, and
we shall have to consider furthermore the cases $\alpha_2 > 1$ and $\alpha_2 < 0$.

Without loss of generality, the two tangent parabolae are given by $y = 0$ and $y = x^2$.
The third parabola is 
$$ y = \alpha_2(x - a_1)(x - a_2),\ a_1,a_2 \in \R^*, a_1 < a_2 \;.
$$

If $\alpha_2 >0$, again using the reflection $x \leftrightarrow -x$, we are lead to only 3
possibilities which are degenerate cases of the preceding ones (by possibly exchanging the
roles of $\alpha$ and $\beta$).  In fact, if $a_1,a_2$ have opposite signs and $\alpha_2
\leq 1$, then the two parabolae $y = x^2$, $y = \alpha_2(x - a_1)(x - a_2)$ do not meet in
real points.

\begin{figure}[htbp]
\begin{center}
   \epsfysize=5.5cm
  \epsfbox{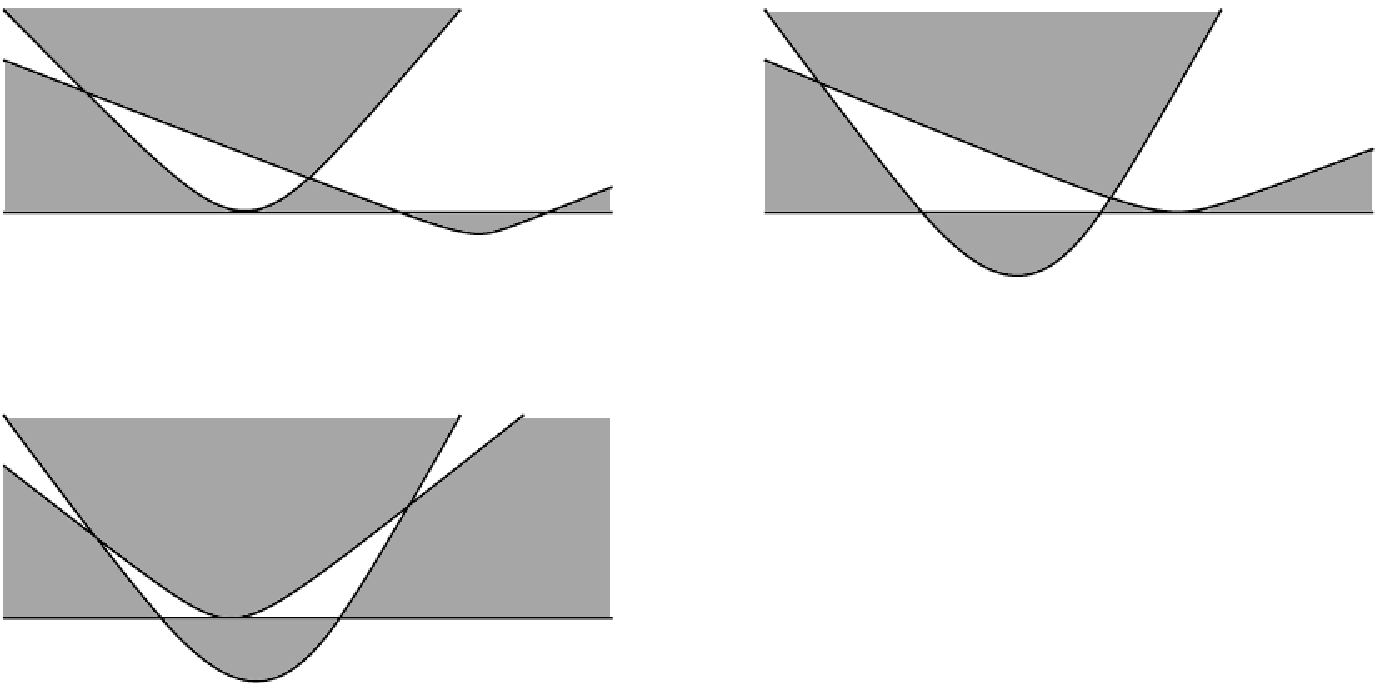}
\caption{Three parabolae, 5 singular points.}\label{fig:3para5points1}
\end{center}
\end{figure}

If $\alpha_2 < 0$ this leads to 2 possibilities, up to reflection $x \leftrightarrow -x$.

\begin{figure}[htbp]
\begin{center}
   \epsfysize=2.3cm
  \epsfbox{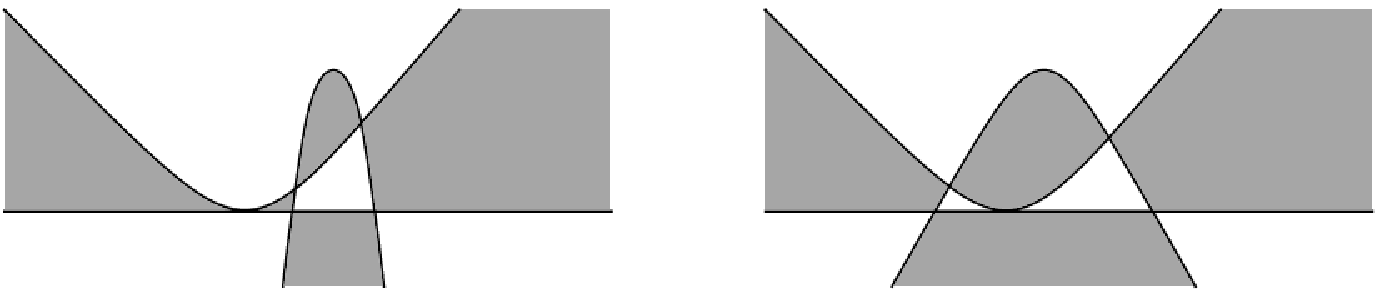}
\caption{Three parabolae, 5 singular points.}\label{fig:3para5points2}
\end{center}
\end{figure}

\end{itemize}

\subsection*{Two components}

\begin{figure}[htbp]
\begin{center}
   \epsfysize=1.4cm
  \epsfbox{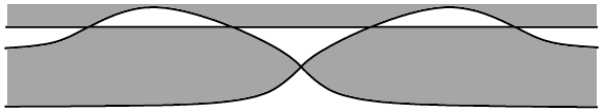}
\caption{The point $O$ is a globally nonseparating $A_1$ singular point.}\label{fig:5points}
\end{center}
\end{figure}

\begin{figure}[htbp]
\begin{center}
  \epsfysize=3.5cm
  \epsfbox{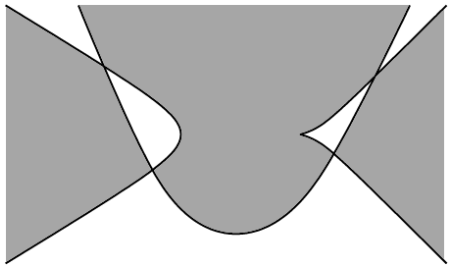}
\caption{The cusp gives rise to a singular point of type $A^-_2$.}\label{fig:cusp}
\end{center}
\end{figure}

Here, we will see that in most cases, any connected component of the topological
normalization of any real double cover ramified over $B$ has at most 4 singular points. There will remain two cases to examine separately, see Figures~\ref{fig:5points} and
\ref{fig:cusp}. Equation~(\ref{eq:B}) becomes 
$$ (y - \alpha(x))(y^2 - \gamma(x)) = 0 \;.
$$

If the bisection $y^2 - \gamma(x) = 0$ is smooth, then  the number of singular points is
bounded from above by 4. Hence we assume the bisection  to have a singular point $O$ at $x =
y  = 0$. To ensure that the bisection and the parabola have 4 real intersection points,  the
polynomial $\alpha(x)^2 -\gamma(x)$ must have 4 distinct roots. These roots are all supposed
to be real and non vanishing in order that $B$ have 5 singular points. The singular point $O$ is either nodal or cuspidal.

If $O$ is an ordinary double point, the bisection is given by 
$ y^2 -x^2h(x) = 0
$ where the quadratic polynomial $h$ is not a square since the bisection is irreducible.
Changing coordinates, we can assume that the parabola is given by $y = 0$ and the bisection
$C$ by $(y + \alpha(x))^2 -x^2h(x) = 0$. Without loss of generality, $\alpha(0) > 0$.

Observe that the leading coefficient of $h$ is non vanishing since the curve $C$ does not
pass through the vertex of $Q$.

The number of real singularities implies that $h$ is not always negative and $h(0) \ne 0$. 
If $h(0) < 0$, then $O$ will be isolated in $B(\R)$ and gives rise to a globally
nonseparating point of the double covering $X'$, in view of the following.

\begin{oss}\label{isolated}
 Let $\pi \colon X'(\R) \to F(\R)$ be a double cover of a smooth connected real
surface $F(\R)$. If $b$ is an isolated point of the real branch curve $B(\R)$, then either
$p = \pi^{-1}(b)$ is an isolated point of $X'(\R)$, or $p$ is a locally separating point of
$X'(\R)$. If however $B(\R)$ has a component $\Gamma$ of dimension 1, then $p$ is globally
nonseparating. 
\end{oss} 

\begin{proof} Take a path connecting $b$ to $\Gamma$.

\end{proof}

\begin{figure}[htbp]
\begin{center}
   \epsfysize=1.4cm
  \epsfbox{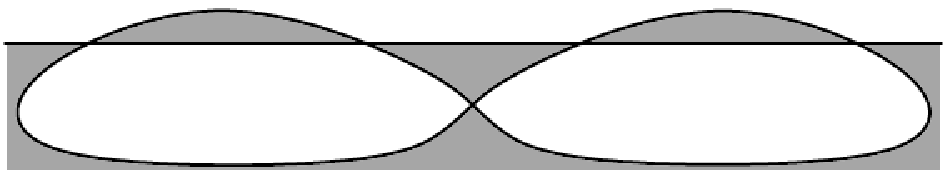}
\caption{Two irreducible components.}\label{fig:4points1}
\end{center}
\end{figure}

If  $h(0) > 0$, and the function $h$ is somewhere negative, observe that $y = 0$ disconnects
the cylinder $\F_2(\R) - \Sigma_\infty(\R)$. Since the polynomial $\alpha(x)^2 - x^2 h(x)$
is assumed to have 4 distinct roots, up to taking a projectivity of $\PP^1(\R)$ sending
$\infty$ to a finite point, we see that there is only one topological possibility, given by
Figure~\ref{fig:4points1}.

\begin{figure}[htbp]
\begin{center}
   \epsfysize=1.4cm
  \epsfbox{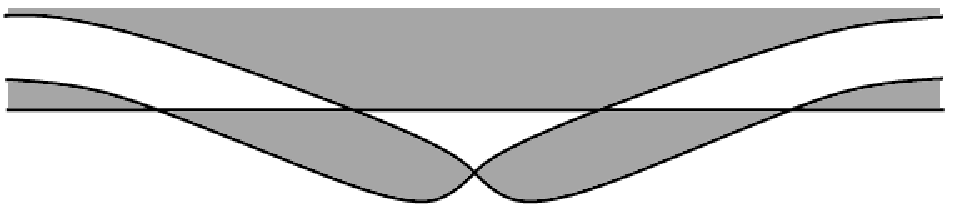}
\caption{Two irreducible components.}\label{fig:4points2}
\end{center}
\end{figure}

If $h(x) > 0$, for all $x$, then $C$ is a double cover of $\PP^1(\R)$, and we can write $C$
as $C^u \cup C^l$, where $C^u$ is the "upper part", $C^l$ the lower part. Because of our
choice $\alpha(0) > 0$, $C^u \cap \{y = 0\} = \emptyset \Rightarrow C \cap \{y = 0\} =
\emptyset$, hence there are two cases: $\#(C^u \cap \{y = 0\} )= 2$, given by
Figure~\ref{fig:4points2}, $\#(C^u \cap \{y = 0\} )= 4$, given by Figure~\ref{fig:5points}.

After we describe the branch curve $B$, observe that we obtain two different surfaces
multiplying the equation of $B$ by $\pm1$.  In Figure~\ref{fig:4points2}, any connected
component of the topological normalization of any double cover will have at most 4 singular
points. In Figure~\ref{fig:5points}, for only one choice of sign, the topological
normalization of the double cover will have a connected component with 5 singular points.
For this double cover however the singular point $O$ turns out to be a globally
nonseparating $A_1$ singular point hence does not belong to $\PPP_{X'}$.

If $O$ is a cusp, the equation of $B$ is
$$ (y - \alpha(x))(y^2 -x^3 l(x)) = 0 \;.
$$

Using a dilatation $y \mapsto \lambda y$ and possibly the usual reflection $x \leftrightarrow
-x$, we may assume $l(0) = 1$ and the equation of the bisection becomes
$$ y^2 -x^3 (1 + ax) = 0\;.
$$

To ensure that the bisection and the parabola have 4 real intersection points,  the equation
$\alpha(x)^2 -x^3 (1 + ax) = 0$ must have 4 distinct roots. Possibly changing the line $x =
\infty$ via a projectivity, we may assume that $a > 0$ and indeed $a = 1$. It is easy then
to see that the only possible configuration is given by Figure~\ref{fig:cusp}. 

Recall that we obtain two different surfaces multiplying the equation of $B$ by $\pm1$. For
only one choice of sign the topological normalization of the double cover will have a
connected component with 5 singular points. For this double cover, however, the point $O$
turns out to be of real type $A^-_2$ which does not belong to
$\PPP_{X'}$.

%\medskip
\subsection*{One component} If the trisection is irreducible, then it has at most 4 singular
points, since $B(\C)$ has genus 4.

\begin{proof}[Proof of Proposition~\ref{prop:dps1case}]

We proceed according to the number of irreducible components of $B$,
recalling that the singular points of $X$ correspond to the singular points of $B$.

If $B$ is irreducible, we have already seen that $B$ has at most 4 singular points.

If instead $B$ has 2 irreducible components, and $B$ has strictly more than 4
singular points, we have seen that $B$ has exactly 5 singular points, and
that the complement $\F_2 (\R) \setminus B (\R)$ has one of the topological
configurations of Figures \ref{fig:4points1}, \ref{fig:4points2},  \ref{fig:5points}, \ref{fig:cusp}.

In the case of Figure \ref{fig:4points1} none of the connected components of the 
complement $\F_2 (\R) \setminus B (\R)$ contains more than 4 points.

The same occurs for the case of Figure \ref{fig:4points2}, while for Figure \ref{fig:5points} there is exactly
one connected component $D$ containing the 5 singular points.
However, in this case the nodal point of the bisection yields a globally nonseparating
singular point of $X'$ for the choice of the positivity region which includes $D$.

Similarly, for the case of  Figure \ref{fig:cusp} there is exactly
one connected component $D$ containing the 5 singular points.
However, in this case the cuspidal point yields a point of type $A^-_2$
which does not belong to
$\PPP_{X'}$.

Assume now that $B$ has 3 irreducible components, and at least 5 singular points.

If there are 6 singular points, the complement $\F_2 (\R) \setminus B (\R)$
 has one of the topological
configurations of Figure \ref{fig:6points}, and none of the connected components of the 
complement $\F_2 (\R) \setminus B (\R)$ contains more than 4 points.

An easy inspection of Figures \ref{fig:3para5points1} and \ref{fig:3para5points1} 
reveals that the same holds also in the remaining cases. 

\end{proof}

\begin{prop}[Koll\'ar]\label{cb} Let $X$ be a real conic bundle with $X$ Du Val. 

Then $m_i \leq 4$, $i
= 1,2,\dots,r$. Moreover, if $m_i = 4$, then $\bar{n} (M_i) \cap \PPP_X$ contains 4 $A_1$ points.
Whereas, if $m_i= 3$, then $\bar{n} (M_i)\cap \PPP_X$ contains at least 2 $A_1$.
\end{prop}
\Proof
The assertion $m_i \leq 4 $ is the last assertion of the proof of cor. 9.8 of
\cite{KoIII}.  But the same argument proves indeed what  we have stated above. 
\qed

\begin{proof}[Proof of Theorem~\ref{teo:main}]

Recall that by \ref{lem:m_i} the numbers $m_1, \dots m_r$ of Du Val singular points 
on the connected components of $\overline { X(\R)}$ which
are not of type  $A^-$ and globally nonseparating is an invariant by extremal birational
contractions. Hence, by Theorem \ref{teo:extremal ray} it suffices to consider the case
where $X$ is either a conic bundle or a Del Pezzo surface.

The case of a conic bundle is settled by Proposition \ref{cb}, and by virtue of
Lemma \ref{lem:reducetodps1} it suffices to consider the case where $X$ is a Du Val
Del Pezzo surface of degree 1.

Now it suffices to apply Proposition \ref{prop:dps1case}.

\end{proof}

%%%%%%%%%%%%%%%%%%%%%%%%%%%%%%%%%%%%%%%%
\section{Real rationally connected Threefolds}\label{sec:3fold}
%%%%%%%%%%%%%%%%%%%%%%%%%%%%%%%%%%%%%%%%

This section is devoted to the proof of Corollary \ref{cor:multiple}.
We first of all introduce the concept of a Werther fibration (cf. \cite{HM05b}),
which allows us to set the integer $k$ on an equal footing in both cases 1) and 2).

Let $S^1\times D^2$ be the {\em solid torus}, where $S^1$ is the unit circle $\{u \in \C
\suchthat \vert u \vert = 1\}$ and $D^2$ is the closed unit disc $\{z \in \C, \vert z \vert
\leq 1\}$.  A {\em Seifert fibration} of the solid torus is a differentiable map of the form
$$ f_{p,q}: S^1\times D^2 \to D^2\;,(u,z)\mapsto u^qz^p\;,
$$  where $p$~and $q$ are natural integers, with $p\neq0$ and
$\gcd(p,q)=1$.  Let~$N$ be a $3$-manifold. A \emph{Seifert fibration} of~$N$ is a
differentiable map~$f$ from $N$ into a differentiable surface~$S$ having the following
property.  Every point~$P\in S$ has a closed neighborhood~$U$ such that the restriction
of~$f$ to~$f^{-1}(U)$ is diffeomorphic to a Seifert fibration of the solid torus.

Let~$A^2$ be the half-open annulus $\{w\in\C\suchthat 1\leq|w|<2\}$.  Let~$P$ be the
differentiable $3$-manifold defined by
$ P=\{((w,z)\in A^2\times\C\suchthat |z|^2=|w|^2-1\}
$. Let~$\omega\colon P\ra A^2$ be the projection defined by~$\omega(w,z)=w$.  It is clear
that~$\omega$ is a differentiable map, that~$\omega$ is a trivial circle bundle over the
interior of~$A^2$, and that~$\omega$ is a diffeomorphism over the boundary of~$A^2$.

\begin{df}\label{def:werther} Let~$g\colon N\ra F$ be a differentiable map from a
$3$-manifold~$N$ without boundary into a differentiable surface~$F$ with boundary. The
map~$g$ is a \emph{Werther fibration} if
\begin{enumerate}
\item the restriction of~$g$ over the interior of~$F$ is a Seifert fibration, and
\item every point~$x$ in the boundary of~$F$ has an open neigborhood~$U$ such that the
restriction of~$g$ to~$g^{-1}(U)$ is diffeomorphic to the restriction of~$\omega$ over an
open neighborhood of~$1$ in~$A^2$.
\end{enumerate}
\end{df}

This definition was introduced in \cite{HM05b}, and is motivated by the following theorem.

\begin{teo}\cite[Theorem 2.6]{HM05b}\label{teo:werther} Let~$N$ be a $3$-dimensional
compact manifold without boundary. Then $N$ is a Seifert fibred manifold or a connected sum
of finitely many lens spaces if and only if  there is a Werther fibration $g\colon N\to F$
over a compact connected differentiable surface~$F$ with boundary. Furthermore $N$ is
Seifert fibred if and only if there exist such a map $g\colon N\to F$ with $\partial F =
\emptyset$.
\end{teo}

Thanks to the Minimal Model Program over $\mathbb{R}$ (\cite{KoII}), the original  setting
for $f \colon W \to X$ in Corollary~\ref{cor:multiple}  is replaced by the following: $W$ is
a real projective $3$-fold with  terminal singularities such that $K_W$ is Cartier along
$W(\mathbb{R})$, $W(\R)$ is a topological 3-manifold, and $f \colon W \to X$ is a rational
curve  fibration over $\mathbb{R}$ such that $-K_W$ is
$f$-ample.

The following result relates the connected components of $W(\R)$ with
the connected components of the topological normalization $\overline { X (\R) }$.

\begin{prop}\cite[Cor.~6.8]{KoIII}\label{prop:norm} Let $W$ be a real projective $3$-fold
with terminal singularities such that $K_W$ is Cartier along $W(\R)$. Let $f \colon W \to X$
be a rational curve fibration over $\R$ such that $-K_W$ is $f$-ample. Let $N \subset W(\R)$
be a connected component. Then $f(N)$ intersects only one of the connected components of
$X(\R) \setminus \operatorname{Sing} X$.
\end{prop}

Let $\overline{n} \colon \overline{W(\R)} \to W(\R)$ be the topological normalization.
The following is the key result which relates the integer $k(N)$ which was
defined above to the numbers $m_i$ of the singularities in $\PPP_X \cap M_i$.

\begin{prop}\cite[Th.~8.1(6)]{KoIII}\label{prop:g} Let $W$ be a real projective $3$-fold
with terminal singularities such that $K_W$ is Cartier along $W(\R)$. Let $f \colon W \to X$
be a rational curve fibration over $\R$ such that $-K_W$ is $f$-ample. Let $N$ be a
connected component of the topological normalization $\overline{W(\R)}$ and assume that $N$
is an orientable topological 3-manifold. Then there exists a small perturbation $g \colon N
\to F$ of $f\vert_{\overline{n}(N)}$ which is a Werther fibration. Furthermore, there is an
injection from the set of multiple fibres of $g$  to the set of singular points of $X$
contained in $f(\overline{n}(N))$ which are of real type $A^+_{\mu}$, $\mu \geq 1$. If
$\partial F = \emptyset$, then $g$ is a Seifert fibration. If $\partial F \ne \emptyset$,
then $N$ is a connected sum of lens spaces and the number of lens spaces is equal to the
number of multiple fibres of $g$.
\end{prop}

\begin{proof}[Proof of Corollary \ref{cor:multiple}]

We want to show that, for each component $N$ of  $\overline{W(\R)}$,
we have   $k(N) \leq 4$. From the above Proposition \ref{prop:g}
it follows that $k(N)$ is the number of multiple fibres of the Werther 
fibration, hence it suffices to bound the number  of singular points of $X$
contained in $f(\overline{n}(N))$ which are of real type $A^+_{\mu}$. If $f(\overline{n}(N))$ is not a connected component of $X(\R) \setminus \operatorname{Sing} X$, then from \cite[8.2]{KoIII}, $N$ is a connected sum of lens spaces and $f(\overline{n}(N))$ may contain some globally nonseparating singular points of type $A^+_1 \sim A^-_1 $. These produce double fibres for $g$, which however correspond to lens space summands $\mathbb{P}^3(\mathbb{R})$. These summands are excluded by the maximality of $a$ in the definition of $k(N)$.
Thus, by Proposition \ref{prop:norm}, it suffices to bound the number  of singular points of $X$
contained in $f(\overline{n}(N))$ which are of real type $A^+_{\mu}$ and globally separating.
Since however these points are a subset of  $\PPP_X \cap M_i$, for some $i \in \{1, \dots r \}$,
the desired inequality follows from Theorem \ref{teo:main}.

\end{proof}

%%%%%%%%%%%%%%%%%%%%%%%%%%%%%%%%%%%%%%%%
\section{Seifert fibrations over a torus}
%%%%%%%%%%%%%%%%%%%%%%%%%%%%%%%%%%%%%%%%

This section is devoted to the proof of Theorem \ref{teo:torus}.

\begin{lem}\label{lem:conic}
Let $r \colon X \to \PP^1_\R$ be a real conic bundle. Suppose that $X$ is a Du Val surface. If $\overline{X(\R)}$ has
a connected component $M$ diffeomorphic to $S^1\times S^1$, then $r$ is smooth along $X(\R)$ and $X(\R) \sim \overline{X(\R)} \sim S^1\times S^1$.
\end{lem}

\begin{proof} We first want to show that $r_M := r\circ \overline{n}\vert_M$ is surjective and  that $\overline{X(\R)} \sim S^1\times S^1$. Assume that $r_M$ is not surjective. Then $\operatorname{Im}(r_M)$ is
homeomorphic to a segment $[a,b]$.  The fibres $r_M^{-1}(a)$ and $r_M^{-1}(b)$ are the ends
of $r_M^{-1}(a,b)$ and they have a (punctured) tubular neighbourhood which is homeomorphic to an
annulus. This shows that $r_M^{-1}(a)$, $r_M^{-1}(b)$ are connected.  The fibre
$r_M^{-1}(a)$ is a simplicial complex of dimension $\leq 1$, and if $S^1 \subset
r_M^{-1}(a)$, then $S^1$ has a (punctured) tubular neighbourhood which is connected, contradicting the
orientability of $M$. Hence $r_M^{-1}(a)$, $r_M^{-1}(b)$ have Euler number 1. 

It suffices to
show that each fibre $r_M^{-1}(t)$, $t \in (a,b)$ has Euler number
$\geq 0$ and we obtain a contradiction to $e(M) = 0$. Looking at the normal forms for singular points of type $A_\mu$ for conic
bundles (given in  \cite[Proof of Cor. 9.8]{KoIII}), we see that every fibre of
$\overline{r} := r \circ \overline{n}$ is either a circle (and then $r$ is smooth on the fibre), or a point, or an interval.  Thus $r_M$ is surjective. 

Again, the Euler number argument shows that all fibres of $\overline{r}$ are circles, hence $r$ is smooth on $X(\R)$ and $M \sim \overline{X(\R)} \sim X(\R)$. 
\end{proof}

\begin{prop}\label{prop:torus-connected} Let $X$ be a real Du Val surface which is rational
over $\C$. Assume that $X(\R)$ contains only singularities of type $A^+_\mu$ and that
$\overline{X(\R)}$ contains a connected component diffeomorphic to $S^1\times S^1$. Then
$\overline{X(\R)}$ is connected, thus  $\overline{X(\R)} \sim S^1\times S^1$. Furthermore,
there is a minimal model of $X$ which is a real conic bundle over $\PP^1$.
\end{prop}

\begin{proof} Let $X$ be as above.  The blow-up of a point of type $A^+_\mu$ for $\mu \geq
2$ induces a homeomorphism between the real parts.  Thus there is a surface $Z$ such that
all singular points are of type $A_1$ and $\overline{Z(\R)}$ is homeomorphic to
$\overline{X(\R)}$.  Let $Z^*$ be a Du Val minimal model of $Z$.  Then by
\cite[Th.~9.6]{KoIII}, $\pi \colon Z \to Z^*$ is the composition of inverses of weighted
blow-ups (of smooth points). Hence $\pi$ is an isomorphism with the exception of a finite
number of smooth points $p_1,\dots ,p_s \in Z^*$ at which one takes the weighted blow-up
which, in local coordinates $(x,y)$ around $p_j$, has the form $\{ xu - vy^2 \} \to
\{(x,y)\}$. Since the weighted blow-up produces globally nonseparating points, there is a
bijection between the connected components of $\overline{Z(\R)}$ and the connected
components $\overline{Z^*(\R)}$.  The weighted blow-up followed by topological normalization
on a disc neighbourhood of $p_j$ has the effect of replacing $p_j$ by a closed segment.
Hence the connected component of $\overline{Z^*(\R)}$ coming from the one of
$\overline{Z(\R)}$ diffeomorphic to $S^1\times S^1$ is again diffeomorphic to $S^1\times
S^1$. Observe again that the singularities of $Z^*$ are only of type $A_1$.

We have two cases:
\begin{enumerate}
\item The minimal model $Z^*$ is a real Del Pezzo surface of degree 1 or 2;
\item  $Z^*$ is a real conic bundle.
\end{enumerate}

In case (1), we have a realization of $Z^*$ as a double cover, and the topological
normalization $\overline {Z^*(\R)}$ can be realized as the real part of a real perturbation
$Z^*_\varepsilon$ of $Z^*$ (by Brusotti's Theorem, [Bru21]).  The surface $Z^*_\varepsilon$ is a
smooth real Del Pezzo surface of degree 1 or 2.  An orientable connected component of such a
surface is a sphere, see e.g. \cite[Chap. 3]{Si89}, so this case does not occur.

Case (2) follows from Lemma~\ref{lem:conic}.
\end{proof}

To prove Theorem~\ref{teo:torus}, we need the conclusion of Proposition~\ref{prop:torus-connected} in a more general setting. First, we give a partial generalisation of Brusotti's theorem in the case of a Du Val Del Pezzo surface. 

\begin{teo}\label{teo:generalizedbrusotti}
Let $X$ be a Du Val Del Pezzo surface. One can obtain, by a global small deformation of $X$, all the possible smoothings of the singular points of $X$.
\end{teo}

\begin{proof}
The main theorem on deformations of compact complex spaces was proven in \cite{grauert}. Good references are \cite{palamodov} and \cite{pinkham}.
The tangent space to  $\operatorname{Def} (X)$ is given by $\operatorname{Ext}^1 (\Omega^1_X, \hol_X)$, see \cite[Cor.~1.1.11]{sernesi}. The obstruction space $\operatorname{Ob} (X)$ is given by $\operatorname{Ext}^2 (\Omega^1_X, \hol_X)$, see \cite[Prop.~2.4.8]{sernesi}.

By the local to global spectral sequence for $\operatorname{Ext}$, we have the following exact sequence
$$
H^1 (X,\Theta_X) \ra  \operatorname{Ext}^1 (\Omega^1_X, \hol_X) \ra  H^0 (\sE xt^1 (\Omega^1_X, \hol_X) ) \ra  H^2 (X,\Theta_X) \ra 0
$$

and  $ \operatorname{Ob} (X) =  H^2 (X,\Theta_X)$.

Therefore the vanishing $H^2 (X,\Theta_X) = 0$ implies:

\begin{itemize}
\item The local deformation space is smooth
\item Global deformations map onto local deformations.
\end{itemize}

We use a calculation by Burns and Wahl \cite[Prop. 1.2]{BW74}, to the effect that, if $S$ is
the minimal resolution of the Du Val singularities of $X$, then
$ p_*  (\Theta_S)=  \Theta_X$.
Whence,  $ H^2 (X,\Theta_X) = H^2 (S,\Theta_S) .$

But the dual space of $H^2 (S,\Theta_S) $ is $ H^0 (S,\Omega^1_S(K_S))$. The conclusion follows from
$ H^0 (S,\Omega^1_S(K_S)) = 0$ since $ H^0 (S,\Omega^1_S) = 0$
and
$ H^0 (S,\hol_S(- K_S)) \neq 0$.

\end{proof}

\begin{lem}\label{lem:localsmoothing}
Consider a real singular point of a surface $X$, of local equation
$z^2 = f (x,y) $ where $f$ vanishes at the origin and has there an
isolated singular point which we assume to be a nonisolated real point.
Then the topological normalization of $X(\R)$ is locally homeomorphic to the real part  $X_{\varepsilon}(\R)$ of the surface $X_{\varepsilon}$ with equation $ z^2 = f (x,y) - \varepsilon $, for $\varepsilon$
sufficiently small and positive.
\end{lem}

\begin{proof}
The real curve $f(x,y) = 0$ has $2m$ arcs entering into the singular point,
ordered counterclockwise, and the region of positivity consists of $m$ sectors,
which alternate themselves to the $m$ sectors of negativity. Furthermore, we have $m \ne 0$ because the origin is a nonisolated real point of the curve. The smooth curve
$f (x,y) =  \varepsilon $ determines  $m$ domains of positivity whose closure
is homeomorphic to the closure of the corresponding sector of positivity of
$f(x,y)$ (where it is contained).
It follows right away that the double cover $ z^2 = f (x,y) - \varepsilon$ replaces
the singular point by $m$ points, one for each connected component of $ X(\R) \setminus \{0\}$.
\end{proof}

\begin{prop}\label{prop:torus-connected-bis} Let $X$ be a real Du Val surface which is rational
over $\C$.  Assume that all locally separating singularities are gobally separating and that
$\overline{X(\R)}$ contains a connected component diffeomorphic to $S^1\times S^1$. Then
$\overline{X(\R)}$ is connected, thus  $\overline{X(\R)} \sim S^1\times S^1$. Furthermore,
there is a minimal model of $X$ which is a real conic bundle over $\PP^1$.
\end{prop}

\begin{proof}  The minimal resolution of a singular point of type $A_{\mu}^- ,\ \mu$ even, induces a homeomorphism between the real parts, thus as in the proof of \ref{prop:torus-connected}, there is a surface $Z$ such that
\begin{itemize}
\item all singular points are of type $A_1$, or of type $A_{\mu}^- , \mu > 1$, $\mu$ odd, or not of type $A$, and 
\item $\overline{Z(\R)}$ is homeomorphic to $\overline{X(\R)}$.  
\end{itemize}
Let $Z^*$ be a Du Val minimal model of $Z$.  Suppose that the minimal model $Z^*$ is a real Del Pezzo surface of degree 1 or 2. We have then a realization of $Z^*$ as a double cover and we can apply Remark~\ref{isolated} to exclude singular points which are isolated real points of the branch curve. By Lemma~\ref{lem:localsmoothing} and Theorem~\ref{teo:generalizedbrusotti}, the topological
normalization $\overline {Z^*(\R)}$ can be realized as the real part of a real perturbation
$Z^*_\varepsilon$ of $Z^*$. The surface $Z^*_\varepsilon$ is a
smooth real Del Pezzo surface of degree 1 or 2.  An orientable connected component of such a
surface is a sphere, so this case does not occur.

Hence $Z^*$  is a real conic bundle and the conclusion follows from Lemma~\ref{lem:conic}.
\end{proof}

\begin{proof}[Proof of Theorem~\ref{teo:torus}]  
The component $N$ of $W(\R)$ is Seifert fibred hence $f(N)$ is the closure of a connected component of $X(\R) \setminus \operatorname{Sing} X$ (see the statement in \cite[8.2]{KoIII} "so we are in the case (4)").
In the proof of [loc. cit. 6.8], using [loc. cit. 4.3],
Kollar claims that  $f(N)$ cannot map to both local components of a locally separating singularity.
Whence, this singularity must be globally separating. Thus all singularities of $f(N)$ which are locally separating are gobally separating. 
We are now in the situation of Proposition~\ref{prop:torus-connected-bis} whence $\overline{X(\R)} \sim S^1\times S^1$. 
Furthermore, the minimal model $Z^*$ is a real conic bundle and Lemma~\ref{lem:conic} gives that $Z^*$ is smooth along $Z^*(\R)$, thus $Z(\R)\cap \PPP_Z = \emptyset$, hence $X(\R) \cap \PPP_X = \emptyset$.

Applying Proposition~\ref{prop:g} on a minimal model of $W \to X$, we get the
conclusion.
\end{proof}

\newpage

\noindent {\bf Author's address:}

\bigskip

\noindent  Fabrizio Catanese\\ Lehrstuhl Mathematik VIII\\ Universit\"at Bayreuth,
NWII\\ D-95440 Bayreuth, Germany

e-mail: Fabrizio.Catanese@uni-bayreuth.de

\noindent Fr\'ed\'eric Mangolte\\ Laboratoire de Math\'ematiques\\ Universit\'e de Savoie\\
73376 Le Bourget du Lac Cedex, France

e-mail: Frederic.Mangolte@univ-savoie.fr

\end{document}